\documentclass[12pt]{article}
\usepackage{amsmath,epsfig,amssymb,float,graphicx,amsfonts, amsbsy}
\usepackage{xspace,colortbl}

\title{Comparing and interpolating distributions on manifold}
 
\author{\Large Nikolay H. Balov
\footnote{Florida State University, Department of Statistics}\\
{\small\texttt{balov@stat.fsu.edu}}}

\begin{document}

\maketitle

\begin{abstract}

We are interested in comparing probability distributions defined on Riemannian manifold. 
The traditional approach to study a distribution relies on locating its mean point and finding 
the dispersion about that point. On a general manifold however, 
even if two distributions are sufficiently concentrated and have unique means, 
a comparison of their covariances is not possible due to the difference in local parametrizations. 
To circumvent the problem we associate a covariance field with each distribution and 
compare them at common points by applying a similarity invariant function on their representing matrices. 
In this way we are able to define distances between distributions. 
We also propose new approach for interpolating discrete distributions 
and derive some criteria that assure consistent results.
Finally, we illustrate with some experimental results on the unit 2-sphere. 

\end{abstract}

%%%%%%%%%%%%%%%%%%%%%%%%%%%%%%%%%%%%%%%%%%%%%%%%%%
{\section{Introduction}}

The problem of comparing distributions defined on a non-Euclidean space or to be more specific,  
a Riemannian manifold, becomes increasingly important. 
A typical example of non-trivial manifold is the unit 2-sphere $\mathbb{S}^2$, 
which is the domain of our experiments in this work.
In this sense, our study has as main application problems from  
directional statistics, a branch of statistics dealing with directions and rotations in $\mathbb{R}^3$.

Pioneers in the field are Fisher, R.A.(1953) and von Mises.
In recent years directional statistics proved to be useful in variety of disciplines like 
shape analysis \cite{mardia-dirs}, geology, crystallography \cite{krieger-crystal}, 
bio-informatics \cite{mardia-protein} and data mining \cite{bahlmann-hand}.
Most of the practitioners in these fields use parametric distributions to model directional data, like 
von Mises-Fisher distribution and Fisher-Bingham-Kent(FBK) distributions.

There are application areas however, where parametric models are insufficient. 
A recent example is provided by medical imaging community. 
In a new technique based on MRI and called High Angular Resolution Diffusion Imaging (HARDI), 
the data is represented by Orientation Distribution Functions (ODFs) which are nothing but 
discrete distributions on the unit 2-sphere. 
These distributions by their nature are multi-modal - they are not concentrated about a particular direction.
They do not follow a parametric model and even if they do the eventual model would be too complicated to be efficient.
Consequently, a non-parametric approach is more natural in processing ODFs. 

In analysis of HARDI data researchers first have to solve the problem of registration between 
different volumes of ODFs, corresponding to the images of different subjects.
For this purpose they need models and algorithms for interpolation between ODFs. 
Second, researchers are interested in comparing different groups of subjects 
using HARDI imaging. Usually, a statistical procedure is employed and hypotheses are tested. 
However, comparison between volumes requires comparison between corresponding ODFs and 
no standard method for this is available.
A third problem in processing HARDI data is building connectivity paths for 
a given volume. Once again we need a consistent way to interpolate between ODFs in order to 
follow an optimal propagating direction.

There are no many choices for interpolation procedure beyond the simplest linear one. 
A recent alternative, using the square root representation of probability mass functions, 
was proposed by Srivastava(2007) and implemented in \cite{chiang-mapping}. 
No existing solution though respects the geometry of the underlying domain. 

In conclusion, we need more models and non-parametric procedures for 
comparing and interpolating distributions on the sphere and on Riemannian manifolds in general.
Approaches that address the non-Euclidean nature of the random variables and 
provide adequate solutions.
It is the main subject of this paper to draw some new directions for searching  
of possible solutions.

What we propose basically is a generalization of the classical concept of covariance of distribution. 
We allow covariance to be defined with respect to any point of distribution domain and 
by doing so we try to workaround the problem of finding the mean point, 
which might not exist or be ambiguous. 
Also, since compact manifolds like $\mathbb{S}^2$ do not admit global parametrizations, 
we pay special attention to use the correct mathematical tool for describing the covariance. 
We not only point out to the well known fact that covariance can be viewed as a bi-linear operator 
and thus defined as a tensor, but specify the exact variance of this tensor. 
It is important to make a distinction between covariance tensor and metric tensor 
on manifold. A central observation in our approach is that at any point of the domain, 
the product of the metric and covariance tensors is a linear operator on the respected tangent space. 
We call it {\it covariance operator}. Collectively they form a field of operators.
Then we introduce instruments, the so called similarity invariant functions, that can be used 
to study properties of these fields and to manipulate them.

After a formal introduction to the concept of covariance operators in section 2, 
we continue in section 3 with motivating examples showing the advantages of the new approach.
We consider a two-sample location problem of the sphere and 
apply several classical non-parametric tests to solve it. 
Test statistics are based on projections defined by covariance operator fields.

In section 4 we consider the problem of interpolation between distributions on the sphere, 
and discuss and compare several alternatives. 
We also show some examples of interpolation between ODFs. The results are encouraging 
in the possibility of developing new applications for processing HARDI.

Although in all our experiments we stay on the unit sphere , 
the theoretical framework still holds on a general Riemannian manifold and 
this is one of its main strengths.

%%%%%%%%%%%%%%%%%%%%%%%%%%%%%%%%%%%%%%%%%%%%%%%%%%
{\section{Covariance fields}}

%%%%%%%%%%%%%%%%%%%%%%%%%%%%%%%%%%%%%%%%%%%%%%%%%%
{\subsection{Random variables on manifold}}

Let M be a Riemannian n-manifold, $q\in M$ and 
let $Exp_q$ be the exponential map at $q$, $Exp_q: M_q \to M$.
If M is complete, then the exponential map $Exp_q$ is defined on the whole tangent space $M_q$.
Throughout this paper for convenience we will assume that M is a complete Riemannian n-manifold, 
although often that is not necessary.

There is a maximal open set $U(q)$ in $M_p$ containing the origin, where $Exp_q$ is a diffeomorphism. 
Then the set $\mathcal{U}(q) = Exp_q(U(q))$ is called maximal normal neighborhood of $q$. 
On this normal neighborhood the exponential map is invertible and let 
$$
Log_q = Exp_q^{-1} : \mathcal{U}(q) \to M_p
$$
be its inverse, the so called log-map. $Log_q$ is diffeomorphism on $\mathcal{U}(q)$.
We adopt the notation $\overrightarrow{qp} = Log_q p$ in analogy to the Euclidean case, 
$M=\mathbb{R}^n$, where $Log_q p = p - q = \overrightarrow{qp}$.

In particular, for $M=\mathbb{S}^n$ the log-map has a closed-form expression 
$$
\overrightarrow{qp} = \frac{cos^{-1}<p, q>}{(1-<p, q>^2)^{1/2}} (p - <p, q> q),
$$
which greatly simplifies metric related operations on the unit sphere.

The Borel sets on M generated by the open sets on M form a $\sigma$-algebra $\mathcal{A}(M)$ on M. 
Any Riemannian manifold has a natural measure $\mathcal{V}$ on $\mathcal{A}(M)$, 
called {\it volume measure}. In local coordinates $x$ it is given by
$$
dV(x) = \sqrt{|G_x|}dx,
$$
where $G_x$ is the matrix representation of the metric tensor, $|G_x|$ is its determinant and 
$dx$ is the Lebesgue measure in $\mathbb{R}^n$. 

\newtheorem{randcomp_ex1}{Example}
\begin{randcomp_ex1}
Consider the two sphere, $\mathbb{S}^2$, parametrized in geographical coordinates $(\theta,\phi)$. 
Then the metric tensor is represented by 
\begin{equation}
G_{(\theta,\phi)} = \left( \begin{array}{cc}
1 & 0 \\
0 & \cos^2(\theta)
\end{array} \right)
\end{equation}
and the volume form is $V(\theta,\phi) = \cos(\theta)d\theta d\phi$.
\end{randcomp_ex1}

A random variable X on M is any measurable function from a probability space 
$(\Omega,\mathcal{B},\mathcal{P})$ to $(M, \mathcal{A}(M), \mathcal{V})$.
The distribution function F of X is defined as 
$$
F(A) = \mathcal{P}(X^{-1}(A)),\\ A\in \mathcal{A}(M).
$$
If $F$ satisfies 
$$
F(A) = \int_A f(p) dV(p), \forall A\in \mathcal{A}(M),
$$
for almost everywhere continuous (w.r.t. $\mathcal{V}$) function $f$, then 
$F$ is said to be absolute continuous (w.r.t. $\mathcal{V}$) and $f$ is its density ({\it pdf}).

%%%%%%%%%%%%%%%%%%%%%%%%%%%%%%%%%%%%%%%%%%%%%%%%%%
{\subsection{Intrinsic and Extrinsic mean and covariance}}

Let $(M,\rho)$ be a metric space. The {\it Fr\'echet mean} set of a distribution F is the 
set of minimizers of $Q(q) = \int \rho^2(q,p)dF(p)$. It was introduced by Frechet (1948).
If M is a Riemannian manifold M with metric structure $g$, then 
the intrinsic mean of $F$, is the Frechet mean of $(M,d_g)$, where $d_g$ is the geodesic distance.
Karcher(1977) considered the intrinsic mean on M and gave conditions for its existence and uniqueness. 
An alternative to intrinsic mean is the extrinsic one, which is obtained by embedding M 
into a higher dimensional Euclidean space.
We point to the influential paper of Bhattacharya R. and Patrangenaru, V. (2005) where the properties 
of extrinsic and intrinsic means and their relation and asymptotic properties are considered in details.

Once a mean point (intrinsic or extrinsic) is specified, the covariance can be defined as usual after fixing a coordinate system about that point.

To compare two distributions one may first look at their intrinsic means. 
If they differ, the distributions differ, otherwise one may compare further their covariances at the common mean point.
This approach however suffers from at least two drawbacks.
First, if the population mean set is large, then the finite sample intrinsic mean will have 
substantial variance. That will diminish the power of any test for equality of means and more importantly, 
will inevitably require comparing covariances at different points.
Second, the intrinsic mean, provided it exists and it is unique, and the covariance 
alone do not specify completely a distribution. 

Thus, if we want to answer the problem 
of comparing distributions, we need a more informative structure that 
completely represents distributions and that is defined in 
coordinates free manner for seamless manipulation and comparison.

%%%%%%%%%%%%%%%%%%%%%%%%%%%%%%%%%%%%%%%%%%%%%%%%%%
{\subsection{Covariance operators}}

Many parametric families of distributions can be defined as functions on linear operators.
Consider for example the standard normal distribution in $\mathbb{R}^n$ with density
$$
f(x) \propto \exp(-\frac{1}{2}||x-\mu||^2),
$$
where $\mu \in \mathbb{R}^n$ is its mean. Since $||x-\mu||^2 = tr((x-\mu)(x-\mu)')$ and 
the matrix $L(x) = (x-\mu)(x-\mu)'$ defines a linear operator 
$$
L(x)(u,v) = u'L(x)v = [u'(x-\mu)][(x-\mu))'v],\\ u,v\in \mathbb{R}^n,
$$
we can express the density by $f(x) \propto h(L(x))$, 
$
h(T) := \exp(-\frac{1}{2} tr(T)).
$

The von Mises-Fisher and FBK distributions \cite{kent-distribution} on the unit 2-sphere give us other such examples.
For example, the latter is given by density
$$
f(\mathbf{x})=\frac{1}{\textrm{c}(\kappa,\beta)}\exp\{\kappa\boldsymbol{\gamma}_{1}\cdot\mathbf{x}+\beta[(\boldsymbol{\gamma}_{2}\cdot\mathbf{x})^{2}-(\boldsymbol{\gamma}_{3}\cdot\mathbf{x})^{2}]\},
$$
where $\gamma_1$, $\gamma_2$ and $\gamma_3$ are three points on $\mathbb{S}^2$ representing orthonormal directions in $\mathbb{R}^3$.
We have $f(x) \propto h(L_1(x), L_2(x), L_3(x))$,
where 
$$
L_i(x) = G(x)(\overrightarrow{x\gamma_i})(\overrightarrow{x\gamma_i})',
$$
are linear operators at the tangent space at $x\in\mathbb{S}^2$
and 
$$
h(T_1, T_2, T_3) = \exp(\kappa\textrm{ }cos(tr(T_1)) + \beta[cos^2(tr(T_{2})) - cos^2(tr(T_{3}))]).
$$
In fact, any $L: x \mapsto L(x)$ in the presented examples is a field of linear operators 
on tangential spaces.
The concept we are going to introduce generalizes the above observations.

We return to a general Riemannain manifold M with metric $G$.
Fix a point $q\in M$. Recall that the metric $G(q)$ is a co-variant 2-tensor at $M_q$, 
while the quantity $(\overrightarrow{qp})(\overrightarrow{qp})'$ is a 
contra-variant 2-tensor at $M_q$. 
The contraction of their tensor product, $G(q)(\overrightarrow{qp})(\overrightarrow{qp})'$,  
is a (1,1)-tensor, which is nothing but a linear operator at $M_q$.

Now the idea becomes clear. For a distribution $F$ on M, 
a linear operator at $M_q$ can be obtained by taking the expectation of 
$G(q)(\overrightarrow{qp})(\overrightarrow{qp})',p\sim F$.

From now on we will use the standard notation $T^2(M_q)$ for co-variant 2-tensors on $M_q$, 
$T_2(M_q)$ for contra-variant 2-tensors on $M_q$ and $T_1^1(M_q)$ for bi-linear operators on $M_q$.

\newtheorem{randcomp_def1}{Definition}
\begin{randcomp_def1}\label{randcomp_def:covfield}
Let $r:\mathbb{R}^+ \to \mathbb{R}^+$ be a continuous function.
Covariance of distribution $F$ on M at point $q\in M$ is defined by 
\begin{equation}\label{randcomp_eq:covfield}
\Sigma(q) = \int_{\mathcal{U}(q)} (\overrightarrow{qp})(\overrightarrow{qp})' r(||\overrightarrow{qp}||)dF(p) 
\end{equation}
and $\Sigma: q \mapsto \Sigma(q) \in T_2(M_q)$ is called covariance field of $F$.
\end{randcomp_def1}
With $r=1$ we obtain the generic covariance field associated with $F$ 
and this is the default choice.

As noted above, $G(q)\Sigma(q)$ is a linear operator on $M_q$, 
which we call {\it covariance operator}.
Hence, $G\Sigma$ is a field of linear operators on M.
With respect to a coordinate system $x$ at q, $G(q)\Sigma(q)$ is represented by a symmetric and 
positive definite matrix $G_x\Sigma_x$, where $G_x$ and $\Sigma_x$ are the representations of $G(q)$ and $\Sigma(q)$.
In other words, $G\Sigma$ is a field of symmetric and positive definite operators on M.

If $v \in M_q$ has components $v_x$ with respect to $x$, we define 
$$
(G(q)\Sigma(q))v := \Sigma_xG_xv_x
$$
and
$$
<v,(G(q)\Sigma(q))v> := v_x'G_x\Sigma_xG_xv_x.
$$
One can check that indeed the last quantity is invariant to coordinate change at q.

It is worth to mention that for a covariance field $\Sigma$ on M, 
$\Sigma^{-1}$ is also symmetric and positive definite and when it is differentiable,
$\Sigma^{-1}$ introduce a new Riemannian metric on M.

If $\Sigma_1$ and $\Sigma_2$ are two covariance fields on M, then $\Sigma_1\Sigma_2^{-1}$ is 
a field of linear operators, i.e. for any $q\in M$, $\Sigma_1(q)\Sigma_2^{-1}(q) \in T_1^1(M_q)$.

On a complete Riemannian manifold, the problem of minimizing the trace of the default covariance field is equivalent to the 
problem of finding the intrinsic mean of $F$, i.e. 
$$
\mu = argmin_{q\in M} \{ \int_{\mathcal{U}(q)} tr(G(q)(\overrightarrow{qp})(\overrightarrow{qp})')dF(p) = \int_{M} d_g^2(q,p)dF(p) \}.
$$

%%%%%%%%%%%%%%%%%%%%%%%%%%%%%%%%%%%%%%%%%%%%%%%%%%
{\subsection{Similarity invariants}}

Let $Sym_n^+$ denote the space of symmetric and positive definite matrices. 
Since this is the representation domain for covariance operators it is of obvious importance for us.
$Sym_n^+$ attracted the attention of many researchers in the recent years due to 
its non-Euclidean nature and consequently, the variety of research opportunities it provides.
For the purposes of Diffusion Tensor Imaging, 
Fletcher, P. T., Joshi, S., (2007)
and Pennec. X., Fillard, P., Ayache, N (2006) proposed 
the use of {\it affine invariant} distance, while 
Arsigny,V., Fillard, P., Pennec X., and Ayache, N. (2007) proposed the so called 
log-Euclidean distance. A good survey of the available distances and estimators in $Sym_n^+$ 
along with new ones is provided by Dryden, I., Koloydenko, A., and Zhou, D., (2008).
We aim a more general treatment of $Sym_n^+$ and instead of dealing with specific 
matrix functions we define a class of invariants. 
What particular member of this class should be used is application problem specific choice. 

Two matrices $A,B\in Sym_n^+$ are said to be similar if
$$
A = X^{-1}BX, \textrm{ for }X\in GL_n.
$$
Matrix representations of linear operators are similar and thus, this fact holds for 
the representations of $G\Sigma$ and $\Sigma_1\Sigma_2^{-1}$.
Next we define an important class of functions that respect similarity. 

\newtheorem{randcomp_def2}[randcomp_def1]{Definition}
\begin{randcomp_def2}\label{randcomp_def:siminvfunc}
A {\it similarity invariant function} on $Sym_n^+$ 
is any continuous bi-variate $h$ that satisfies
\begin{enumerate}
\item[(i)] $h(AXA',AYA') = h(X,Y)$, $\forall X,Y\in Sym_n^+$ and $A\in GL_n$.
\end{enumerate}
It is a non-negative with a unique root if
\begin{enumerate}
\item[(ii)] $h(X,Y) \ge 0$, $\forall X,Y\in Sym_n^+ \textrm{ and } h(X,Y) = 0 \iff X = Y$.
\end{enumerate}
Moreover, $h$ is called similarity invariant distance, if in addition to (i) and (ii) also satisfies
\begin{enumerate}
\item[(iii)] $h(X,Y) + h(Y,Z) \ge h(X,Z)$, $\forall X,Y,Z\in Sym_n^+$.
\end{enumerate}
\end{randcomp_def2}
Below we list several examples of similarity invariant function we use in our experiments. 

\begin{enumerate}

\item{For a fixed $Z\in Sym_n^+$, the similarity invariant 
$$
h_{trdif}(X, Y; Z) = |(tr(Z^{-1}X - Z^{-1}Y)|,
$$
satisfies (iii) but not (ii). Default choice will be $Z=G^{-1}$, the inverse of the metric tensor representation.
}
\item{The second one is sometimes referred as {\it affine-invariant distance} in $Sym_2^+$, 
see for example 
\cite{ohara-suda-amari}, \cite{forstner-moonen}, \cite{batchelor-tensors}, \cite{fletcher-tensors} and \cite{pennec-tensors}, and it is defined by 
$$
h_{trln2}(X, Y) = \{tr(ln^2(XY^{-1}))\}^{1/2}, X,Y\in Sym_2^+.
$$
Actually, $h_{trln2}$ is not a unique choice for a distance in $Sym_2^+$.
}
\item{
Log-likelihood function gives us another choice for h, 
$$
h_{lik}(X, Y) = tr(XY^{-1}) - ln|XY^{-1}| - n.
$$
It satisfies (i) and (ii) but it fails to satisfy the triangular inequality.
}
\item{
Another interesting choice for $h$ is
$$
h_{lnpr}(X, Y) = \{ln(tr(XY^{-1})tr(YX^{-1}))\}^{1/2},
$$
that satisfies (iii) and 'almost' satisfies (ii):  
$h_{lnpr}(X, Y) = 0 \iff X = cY$, for $c>0$.
}

\end{enumerate}

The concept of covariance fields 
can be used for measuring the difference between distributions on M.
Let $f$ and $g$ be two densities on M and $\Sigma[f]$ and $\Sigma[g]$ be 
their respected covariance fields. 

For a non-negative $h\in \mathcal{SIM}(n)$ we define 
\begin{equation}\label{randcomp_eq:compare_densities_dh}
d_h(f, g) := \int_{M} h(\Sigma[f](p), \Sigma[g](p)) dV(p).
\end{equation}
When M is a compact, the above integral is well defined and finite.
Moreover, if $h(X,Y)$ is a distance function on $Sym_n^+$, 
then $d_h$ will be a distance in the space of densities on M.

Equation (\ref{randcomp_eq:compare_densities_dh}) gives a very general but impractical way to compare 
distributions due to the fact that the integration domain is the whole manifold. 
For application purposes however, one may restrict to a smaller domain 
or perform the comparison on discrete set of points which are of particular interest. 

%%%%%%%%%%%%%%%%%%%%%%%%%%%%%%%%%%%%%%%%%%%%%%%%%%%
{\section{Two-sample location problem on $\mathbb{S}^2$}}

In this section we make an application of covariance operators to non-parametric 
distribution comparison. It will serve more illustration purposes rather than strong application ones.  
The goal is to provide motivating examples showing the new opportunities provided by 
the proposed covariance structure.
We choose to apply simple procedures, as Wilcoxon signed rank and rank sum tests, 
in order to have a good look and intuition of what happens.

Let $\{p_{i,1}\}_{i=1}^m$ and $\{p_{i,2}\}_{i=1}^m$, be random samples from 
distributions $F_1$ and $F_2$ on $\mathbb{S}^2$, respectively, and the two samples be independent of each other.
 
Fix a point $q\in M$ and define 
$$
\eta_i^1 = G(q)(\overrightarrow{qp_{i,1}})(\overrightarrow{qp_{i,1}})', \\ 
\eta_i^2 = G(q)(\overrightarrow{qp_{i,2}})(\overrightarrow{qp_{i,2}})'.
$$
Using tensor notation we can write $\eta_i^l \in T_1^1(\mathbb{S}_q^2$).
The respective sample covariance operators at $q$ are 
$$
\hat L^1(q) = \frac{1}{m} \sum_{i=1}^m \eta_i^1 \textrm{ and } 
\hat L^2(q) = \frac{1}{m} \sum_{i=1}^m \eta_i^2.
$$
We call $q$ {\it observation} point and basically, what we are going to show is how its choice influences 
the inference about $F_1$ and $F_2$. 

Fix a tangent vector $v\in \mathbb{S}_q^2$ and consider following (ordinary) random variables 
$$
\xi_i^1(v) = <v, \eta_i^1(v)>_q \textrm{ and } \xi_i^2(v) = <v, \eta_i^2(v)>_q,
$$
where $<.,.>_q$ is the dot product in the tangent space $\mathbb{S}_q^2$.

\newtheorem{randcomp_def3}[randcomp_def1]{Definition}
\begin{randcomp_def3}\label{randcomp_def:equal_distr}
We say that $F_1$ and $F_2$ have the same location w.r.t. $q\in \mathbb{S}^2$ 
and write $F_1 \cong_q F_2$ if for any $v\in\mathbb{S}_q^2\cong\mathbb{R}^2$, random variables 
$$
\xi^l(v) = <v, (G(q)(\overrightarrow{qX_l})(\overrightarrow{qX_l})')(v)>_q,
\textrm{ }X_l \sim F_l,\textrm{  }l=1,2
$$
have the same median.
\end{randcomp_def3}

Under the hypothesis $H_0: F_1 \cong_q F_2$,  
for any $v$, $\xi_i^1(v)$ and $\xi_i^2(v)$ are random samples from distributions with equal median. 

To test $H_0$ we propose two procedures based on 
the Wilcoxon signed rank test, \cite{hollander-nonparam}, page 36. 
Let $T_{xi}(v)$\footnote{$T_{xi} = \sum_i r_i s_i$, where for $z_i = \xi_i^1-\xi_i^2$, 
$s_i = 1_{\{z_i>0\}}$ and $r_i$ are the ranks of $|z_i|$.} 
be the signed rank statistics based on $\xi_i(v)$'s and 
$T_{xi} = \max\{T_{xi}(v), v\in \mathbb{R}^2\}$.
Then we reject $H_0$ when $T_{xi}$ is sufficiently large.

The second test is based on $T_d$, the Wilcoxon signed rank statistics for the distances 
$$
d_i^1 = tr(\eta_i^1(v)) = d^2(q, p_{i,1}) \textrm{ and } d_i^2 = tr(\eta_i^2(v)) = d^2(q, p_{i,1}),
$$
where $d^2$ stands for the spherical distance, $d^2(q, p) = cos^{-1}(<q,p>)$.

If we choose an orthonormal basis $\{v_s\}_{s=1}^2$ of $\mathbb{S}_q^2 \cong \mathbb{R}^2$ and define 
$$
\xi_{i,s}^1 = <v_s, \eta_i^1 v_s >_q \textrm{ and } \xi_{i,s}^2 = <v_s, \eta_i^2 v_s>_q,
$$
then the following holds: for any $i=1,...,m$ and $l=1,2$
\begin{equation}\label{randcomp_d_sum_xi}
\sum_{s=1}^2 \xi_{i,s}^l = \sum_{s=1}^2 <v_s,\overrightarrow{qp_{i,l}}>^2 = ||\overrightarrow{qp_{i,l}}||^2 = tr(\eta_i^l(v)) = d_i^l.
\end{equation}

It is clear that if $\xi^1(v_1)$ and $\xi^2(v_1)$ have the same median and 
$\xi^1(v_2)$ and $\xi^2(v_2)$ also have the same median, then 
$\xi^1(v)$ and $\xi^2(v)$ will have the same median for every $v$. 

$\{d_i^1\}_i$ and $\{d_i^2\}_i$ are samples from marginals of $F_1$ and $F_2$, 
which under the null hypothesis have the same location. 
As distances, they are invariant to rotation of the samples $p_{i,l}$ on the sphere. 
On the other hand, for any $l$, 
$\{\xi_{i,1}^l\}_i$ and $\{\xi_{i,2}^l\}_i$ follow two marginal distributions 
that can be considered projections of $F_l$ onto two orthogonal axes. As such 
they form more discriminating set of variables than $\{d_i^l\}_i$. 

These observations motivate the following procedure for testing $H_0$.
\newtheorem{randcomp_alg1}{Test Procedure}
\begin{randcomp_alg1}
Let $\{p_{i,l}\}_{i=1}^m$, l=1,2 be two random samples, independent of each other.
\begin{enumerate}
\item Find the operators $\eta_i^1$ and $\eta_i^2$ and set 
$$
\hat L(q) = \hat L^1(q) - \hat L^2(q) = \frac{1}{m} \sum_{i=1}^m \eta_i^1 - \frac{1}{m} \sum_{i=1}^m \eta_i^2.
$$
\item Let $\lambda_s$ and $v_s$ be the eigenvalues and eigenvectors of $\hat L(q)$. \\
Set $\xi_{i,s}^l = <v_s, \eta_i^l v_s >_q$.
\item Calculate statistics $T_{xi, s}$ based on $\xi_{i,s}^l$ and set 
$$
T_{xi} = \max\{T_{xi, 1},T_{xi, 2}\}.
$$
\item Choose a significance level $\alpha$. 
If $pval(T_{xi}) < \alpha/2$ \footnote{We apply Bonferroni correction for the p-value.}, 
reject $H_0$.
\end{enumerate}
\end{randcomp_alg1}

We also employ the rank sum test (Wilcoxon, Mann and Whitney), \cite{hollander-nonparam}, page 106, 
to compare the performance of $\xi_{i}^l$ and $d_{i,l}$ random variables. For the statistics 
$W_{xi}$
\footnote{$W_{xi} = \sum_i r_i $, where for $r_i$ are the ranks of $\xi_{i}^1$ in the joint sample $\{\xi_{i}^1,\xi_{i}^2\}$.}  
and $W_d$, we calculate corresponding p-values using large sample approximation. 
The second test procedure is the same as the first one but uses $W$ instead of $T$ statistics.

Note that if $F_2$ distribution is a rotation of $F_1$ about $q$, then the type II error 
of $T_d$ statistics will be 1, i.e. the power will be 0.

The way of choosing the basic vectors $v_s$ of $\mathbb{S}^2_q$ resembles  
the Principal Component Analysis (PCA) of the operator $\hat L(q)$ and its derivatives like 
Principal Geodesic Analysis (PGA), introduced by Fletcher, 2004.
In the standard setup, PCA is applied on the covariance defined 
at the (extrinsic or intrinsic) mean point. 
However, not only the existence of a mean is not guaranteed, but 
its properties may not be optimal in the context of the test statistic. 
In contrast, in our approach we allow freedom of choosing the observation point $q$ according 
to a criteria favoring that statistic.

\begin{figure}
\centering
\begin{tabular}{cc}
{\includegraphics[scale=0.4]{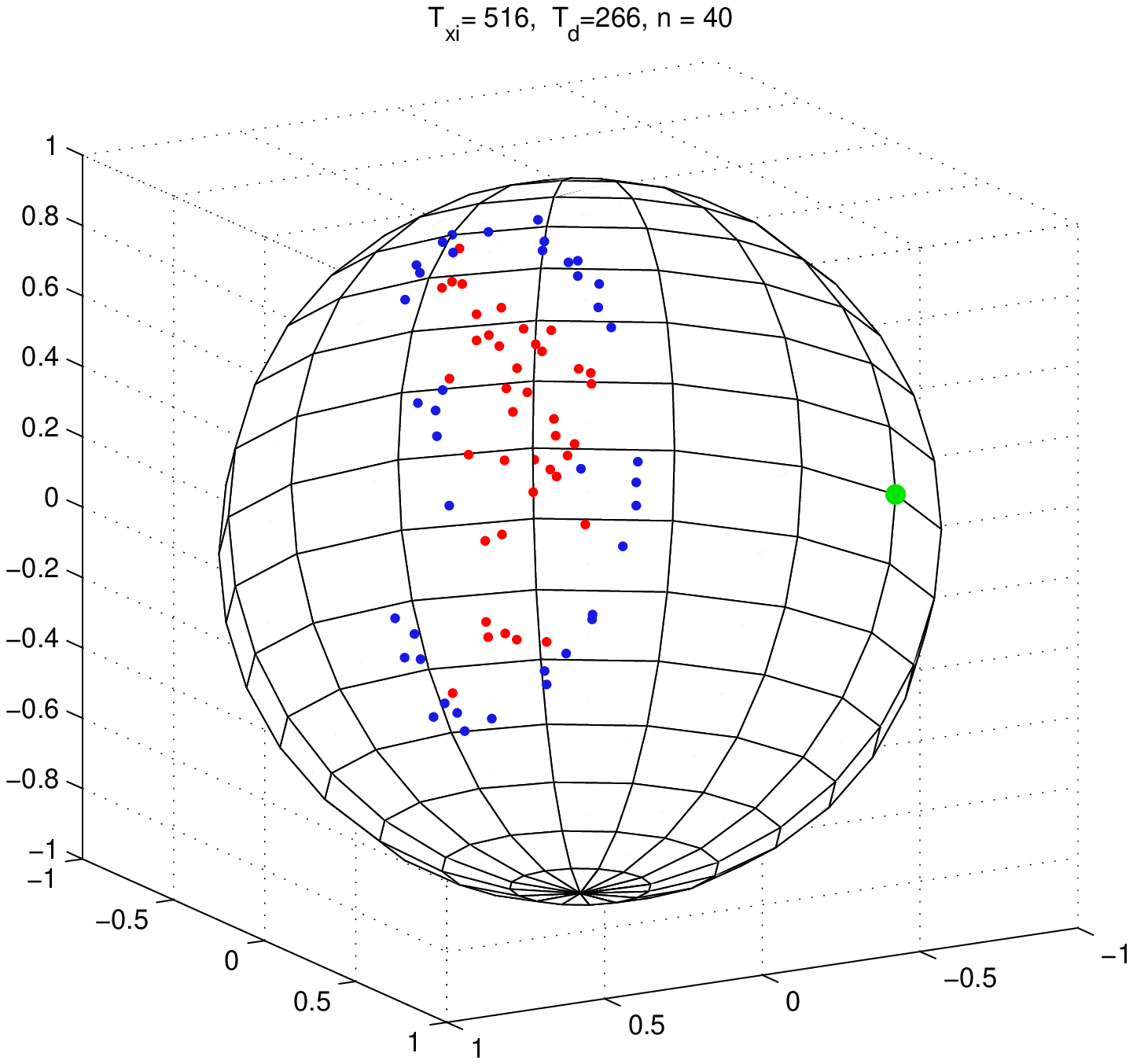}}
{\includegraphics[scale=0.4]{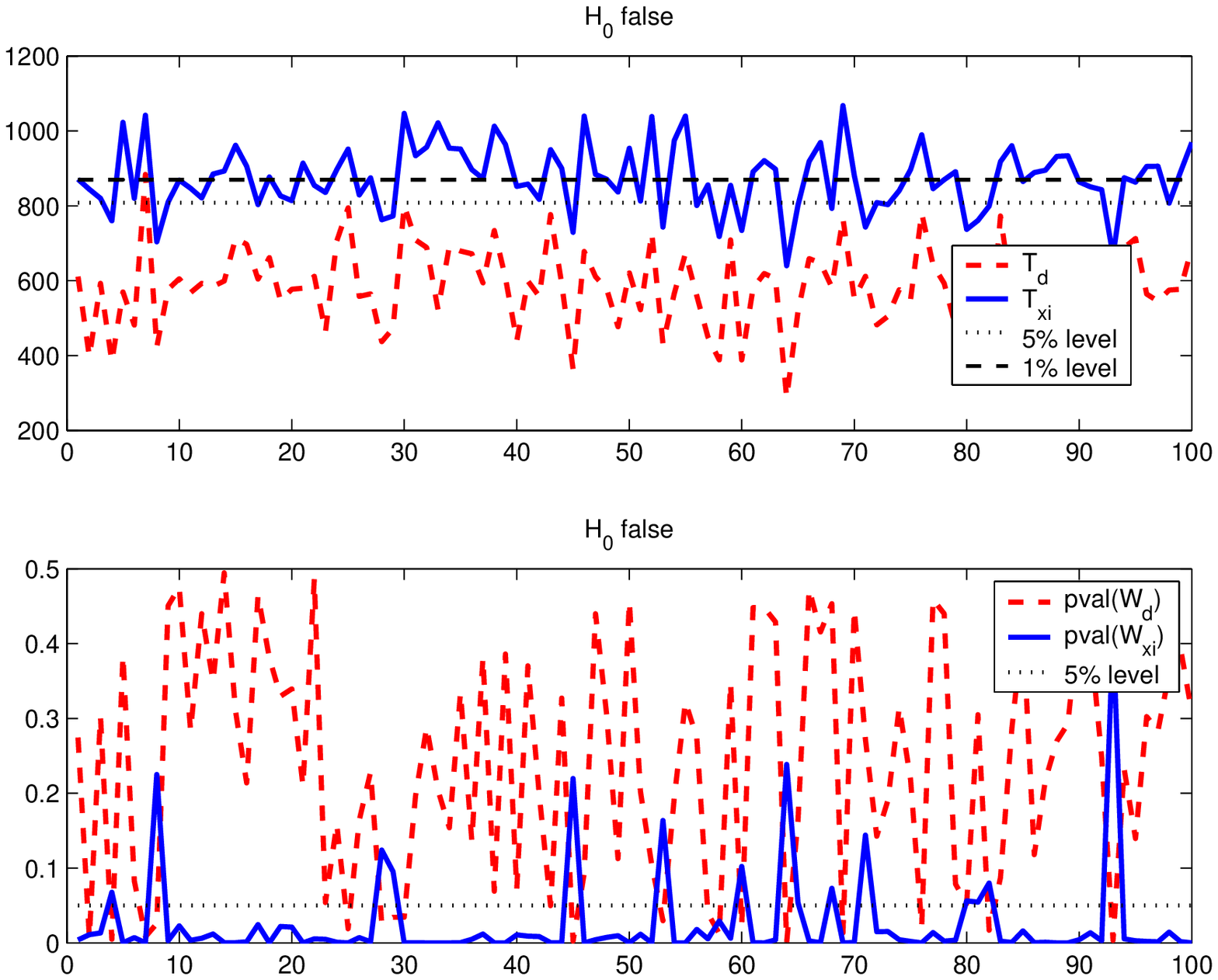}}
\end{tabular}
\caption{
Testing $H_0:F_1 \cong_q F_2$ when $H_0$ is false. 
Sample examples from $F_1$ (red) and $F_2$(blue) are given in the left.
Observation point $q$, shown in green, is fixed and equals distribution parameter $\mu$.
Top right plot shows $T_{xi}$ and $T_d$ statistics along with $1$ and $5$ percentile lines, while 
the bottom right plot shows $W_{xi}$ and $W_d$ by their p-values. 
Test procedure 1 is run 100 times with sample size of 50.
$T_{xi}$ clearly outperforms $T_d$ statistics by rejecting the null hypothesis most of the time and 
the same is true for $W_{xi}$ versus $W_d$.
}
\label{randcomp_fig:1}
\end{figure}
\begin{figure}
\centering
\begin{tabular}{cc}
{\includegraphics[scale=0.4]{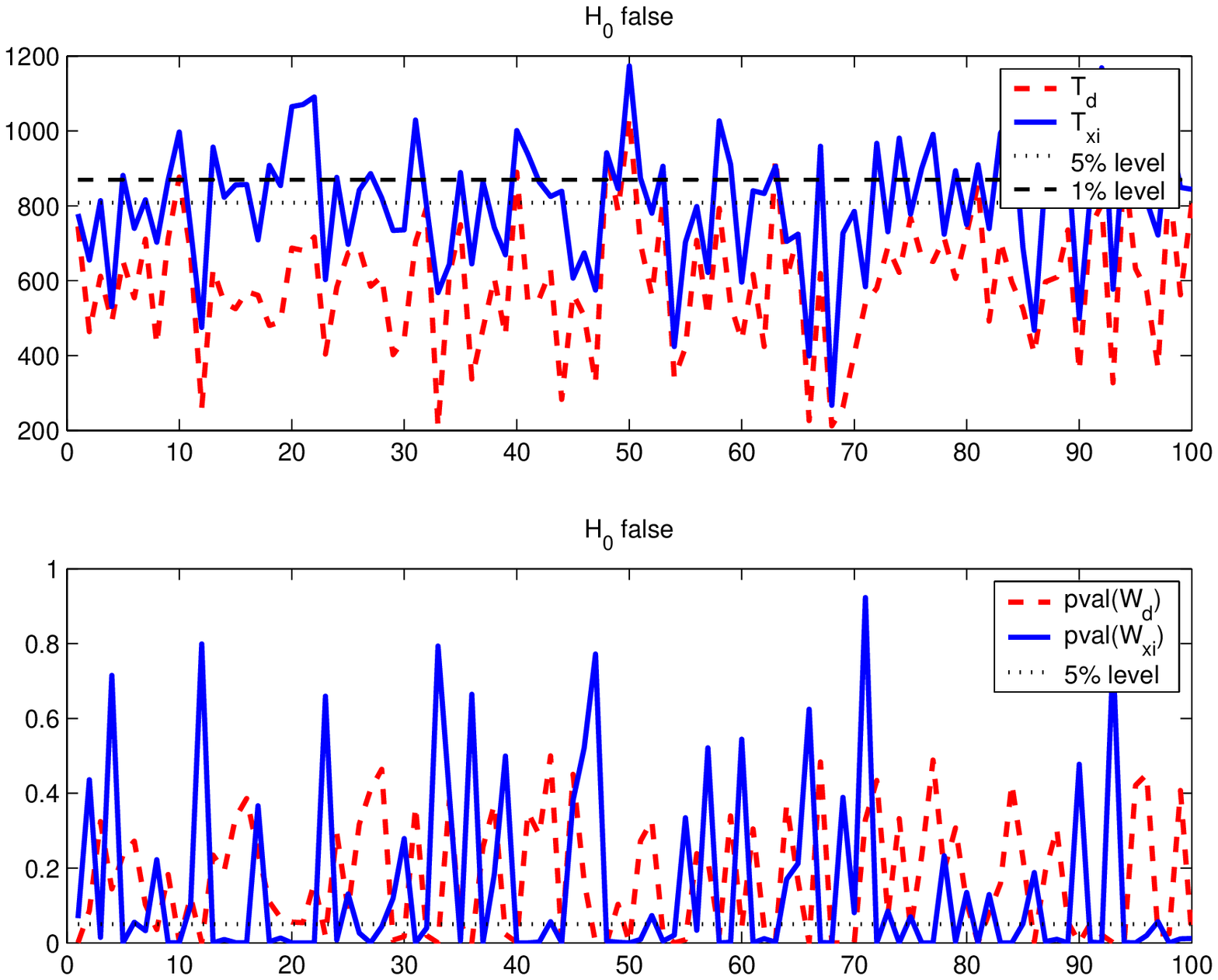}}
{\includegraphics[scale=0.4]{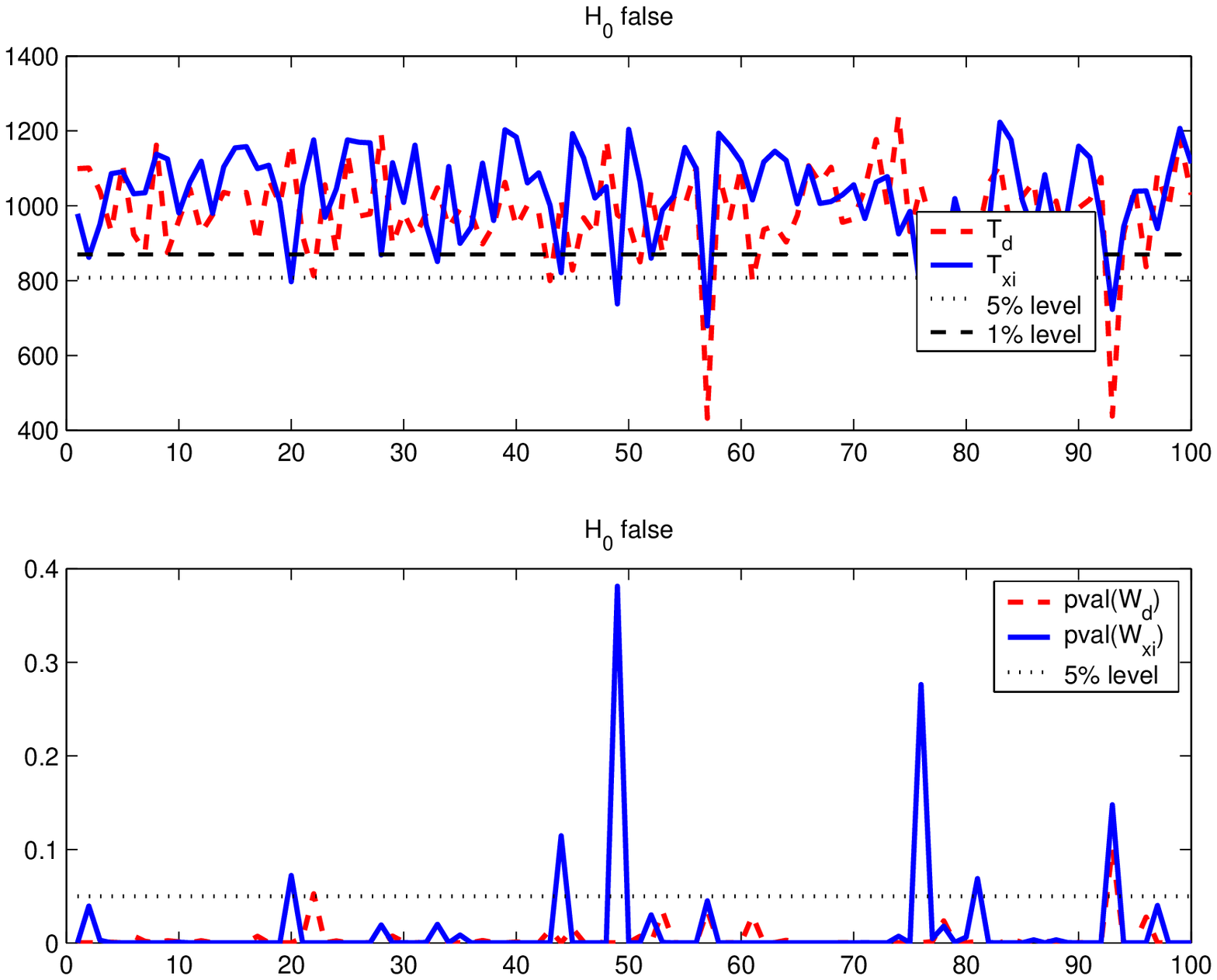}}
\end{tabular}
\caption{Performance of $T$ and $W$ statistics when the observation point $q$ varies.
In the left plot $q$ is chosen uniformly on the sphere. 
The experiment confirms a clear advantage for $T_{xi}$. 
In the right plot preference is given to those observation points 
for which $tr^2(\hat L)$ is large. Now $T_{d}$ is on a par to $T_{xi}$ with both being very high.
Corresponding $W_d$ and $W_{xi}$ also have very significant p-values.
}
\label{randcomp_fig:2}
\end{figure}

Figures \ref{randcomp_fig:1} and \ref{randcomp_fig:2} show some experimental results using 
the proposed procedures for testing $H_0: F_1 \cong_q F_2$.
We consider a family of distributions given by density 
\begin{equation}\label{randcomp_eq:pdf}
f(p; a, \mu) \propto \exp(-(tr(G(\mu)(\overrightarrow{\mu p})(\overrightarrow{ \mu p}))^2 - a)^ 2)
\end{equation}
where $\mu$ is a fixed point (not to be mistaken as a mean) and $a$ is a parameter.
Top plots show Wilcoxon sign rank statistics $T$, while bottom plots show rank sum statistics $W$.
As we see in figure \ref{randcomp_fig:2} left, where the observation point $q$ varies uniformly on $\mathbb{S}^2$, 
for the majority of positions, $T_{xi}$ and $W_{xi}$ achieve higher p-values than $T_d$ and $W_d$.
This result is not isolated and can be repeated for a great variety of distributions 
besides (\ref{randcomp_eq:pdf}).

How the choice of observation point $q$ affects the relative performance of $T$ and $W$ statistics? 
We have that 
$$
d_i^1 - d_i^2 = \sum_{s=1}^2 (\xi_{i,s}^1 - \xi_{i,s}^2) \textrm{ and } \frac{1}{m} \sum_{i=1}^m (\xi_{i,s}^1 - \xi_{i,s}^2) = \lambda_s ,
$$
thus
$$
\frac{1}{m} \sum_{i=1}^m (d_i^1 - d_i^2) = \sum_s \lambda_s.
$$
Therefore, if all $\lambda_s$ are of equal sign, the absolute value of 
the sample expectation of $(d_i^1 - d_i^2)$ will be 
higher than that of $(\xi_{i,s}^1 - \xi_{i,s}^2)$, for all $s$. 
In case when the eigenvalues are of different signs the reverse is expected, 
the absolute value of the sample expectation of $(\xi_{i,s}^1 - \xi_{i,s}^2)$ 
for the maximal $|\lambda_s|$ will be higher than that of $(d_i^1 - d_i^2)$, 
which means that $T_{xi}$ is expected to be higher than $T_d$.
Of course these considerations are only approximate because the tests for $T$ and $W$ statistics 
are based on assumptions on the medians not on the means. Nevertheless, we may take the above as a general observation that 
can be made more formal and rigorous using other appropriate statistical tests. 

\begin{figure}
\centering
\begin{tabular}{cc}
{\includegraphics[scale=0.4]{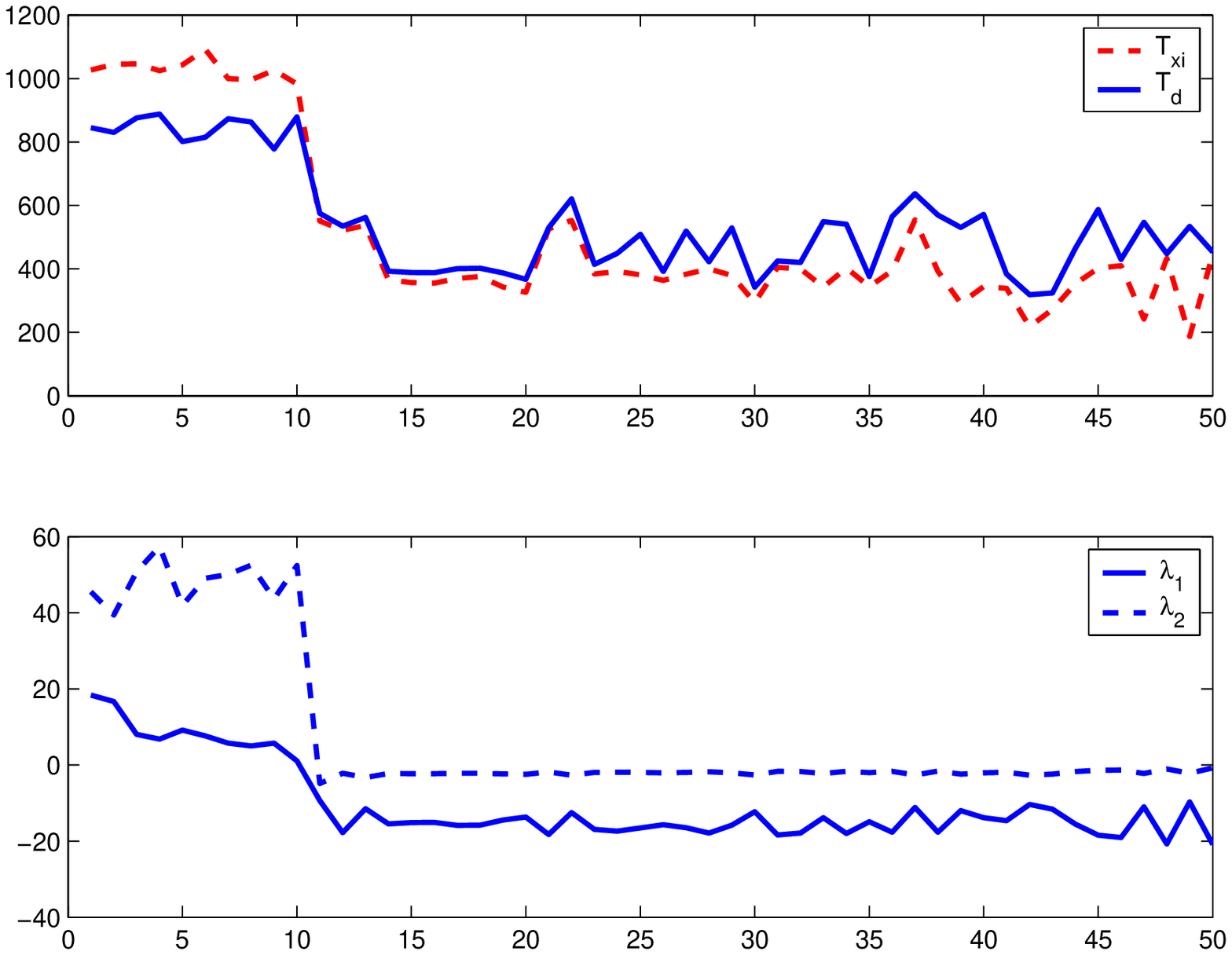}}
{\includegraphics[scale=0.4]{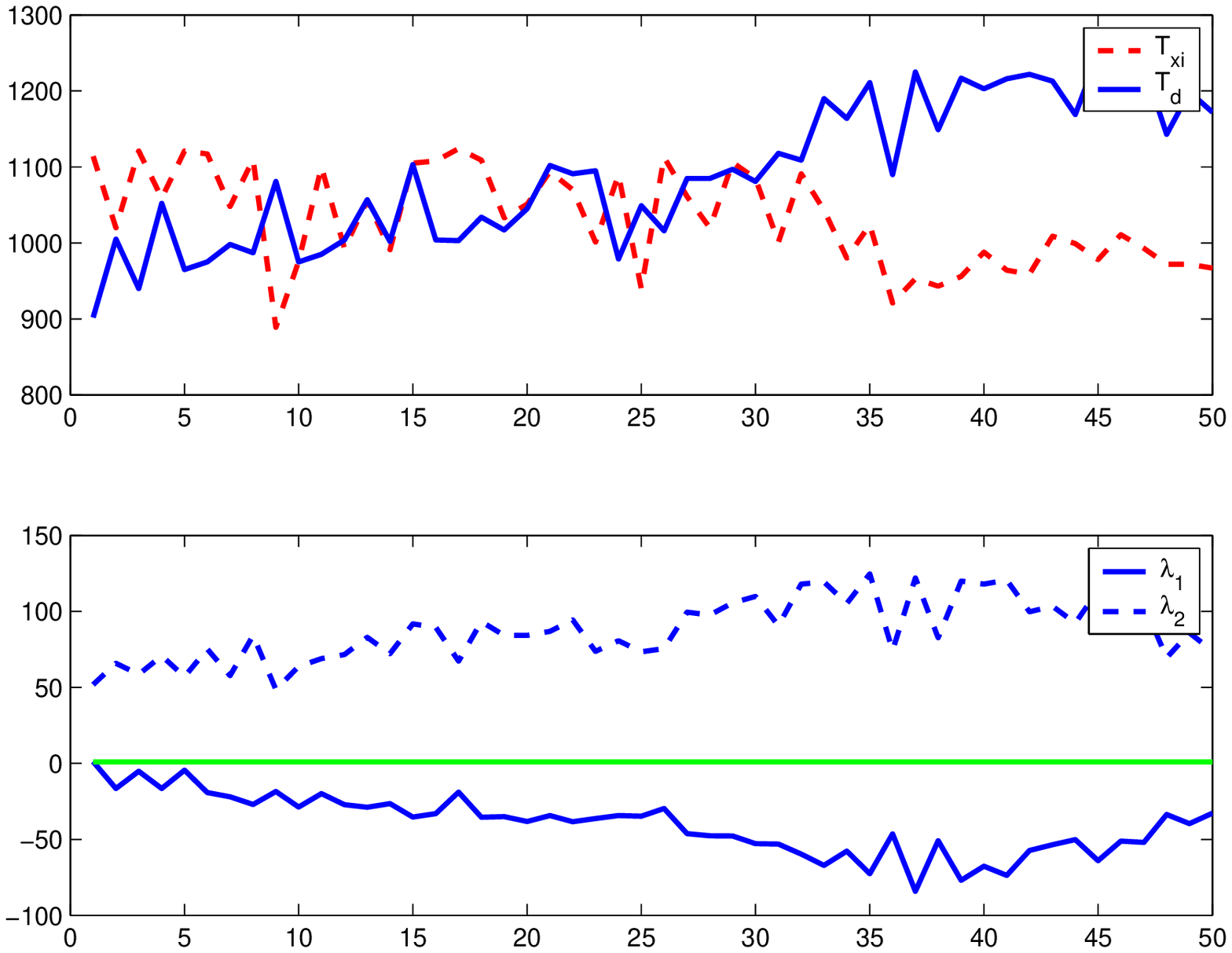}}
\end{tabular}
\caption{Comparing the performances of $T_{xi}$ and $T_d$ statistics (top) 
for a fixed pair of samples by varying the observation point. 
Samples are drawn from (\ref{randcomp_eq:pdf}).
Observation points are ordered decreasingly in $det(\hat L(q))$ in the left plot and 
decreasingly in $tr^2(\hat L(q))$ in the right. 
The bottom plot shows the eigenvalues of $\hat L(q)$. 
We note that $T_d$ is the larger statistics and thus, has lower p-values, only when 
both eigenvalues are strictly positive.
}
\label{randcomp_fig:3}
\end{figure}

We provide some experimental evidence confirming the above expectations. 
For comparison the performance of $T_{xi}$ and $T_{d}$ we use $det(\hat L(q))$. 
We expect for $T_d$ to benefit from positive values of $det(\hat L(q))$ and 
indeed this is the case as seen in figure \ref{randcomp_fig:3} left.
There, for a fixed pair of samples, 
we calculate and compare $T_{xi}$ and $T_d$ statistics at 50 observation points on the sphere. 
Then we sort the results such that $det(\hat L(q))$ decreases. 
In the far left, both $\lambda_1$ and $\lambda_2$ are positive, which leads to a clear advantage for $T_d$. 
Once the sign of $det(\hat L(q))$ goes negative, the situation reverses.

We also expect that at observation points with high values of $tr^2(\hat L)$, 
all statistics to be strong in rejecting a false null hypothesis. 
Some evidence confirming this is shown in figures \ref{randcomp_fig:2} right and \ref{randcomp_fig:3} right.
$tr^2(\hat L)$ is probably the simplest statistics that measures the difference between the two samples and 
it is in fact, an application of the similarity invariant function $h_{trdif}$ as defined in section 2.2.

One can show that $\hat L$ is a continuous field of linear operators on $\mathbb{S}^2$ 
(the proof is beyond the scope of the paper). Therefore, 
if there exists a point $q$ with $det(\hat L(q)) < 0$, 
then that sign is negative on non-vanishing area. 
Only when samples $p_{i,l}$ collectively are highly concentrated, the area $S_+$ where $det(\hat L(q)) > 0$ 
will dominate over $S_-$, the area where $det(\hat L(q)) < 0$. 
In case when $H_0$ is false, we expect that $S_+ < S_-$.

\begin{figure}
\centering
{\includegraphics[scale=0.75]{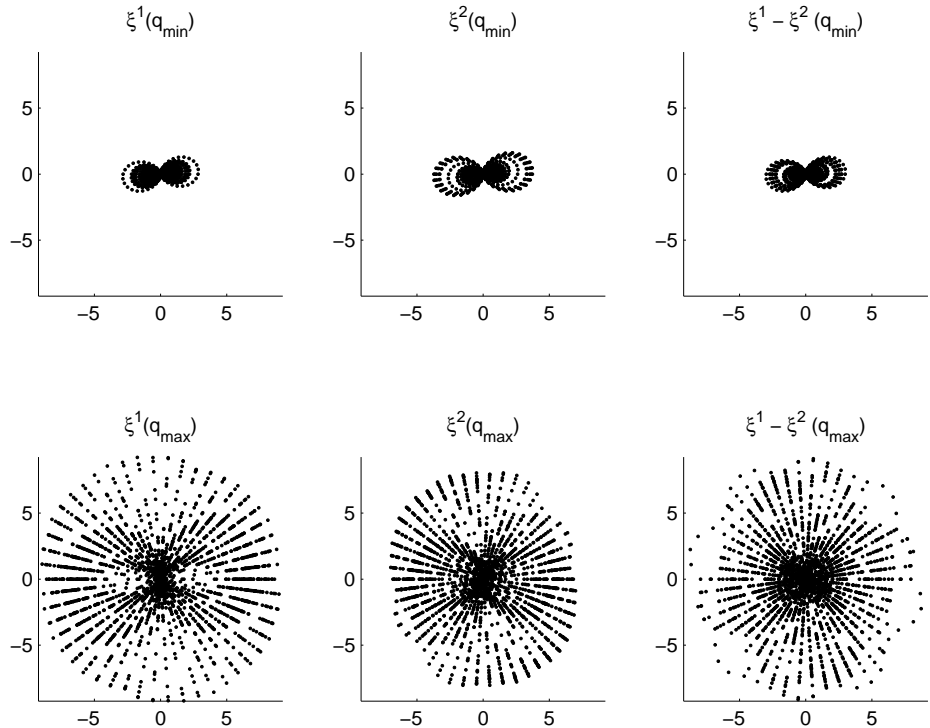}}
\caption{Two samples of size 50 are drawn from $f(p;0.2)$ and $f(p;0.3)$ as given by (\ref{randcomp_eq:pdf}).
Points $q_{min}$ and $q_{max}$ are chosen to minimize and maximize $tr^2(\hat L_1 - \hat L_2)$.
Shown are sample profiles and their difference (right) at these points, 
defined by the projections $\xi = <v, \eta(v)>$ along 50 directions $v$ spanning uniformly $[0,2\pi]$.
Profiles of $\xi^1$ and $\xi^2$ are concentrated and look similar at $q_{min}$, but 
are diffused and very different at $q_{max}$.
These plots visualize clearly the difference between two extreme observation points.
%In the latter case the sample covariance operators $\hat L^1$ and $\hat L^2$ have very different principal directions.
}
\label{randcomp_fig:4}
\end{figure}

Figure \ref{randcomp_fig:4} gives another useful way to visualize 
the sample operator $\hat L$ at different observation points. 
By choosing a point $q$, one can draw the projections $<v,\eta_i^l(v)>_q$ 
for a set of directions $v$ spanning a circle to obtain the so called sample profile. 

In conclusion, choosing an observation point for comparing locations of two distributions 
is an important issue since not all positions provide same test performance.
Position optimality depends on the statistic applied on the covariance operator.
For the projection based statistics we used as examples, optimal observation points can be chosen 
by maximizing the squared trace of the difference of the covariance operators.
 
We also showed that distance based statistics have limited performance 
and in general, employing the whole covariance structure is beneficial. 

We also note that most of the presented results do not depend on the specific geometry of the unit sphere 
and still hold on a general Riemannian manifold.

%%%%%%%%%%%%%%%%%%%%%%%%%%%%%%%%%%%%%%%%%%%%%%%%%%
{\section{Interpolation of discrete distributions on $\mathbb{S}^2$}}

The second application of the covariance operators we are going to consider is 
interpolation between discrete distributions on the unit sphere. 
We suppose that the distributions are defined on a common domain - a fixed set of points on the sphere. 
The approach we propose is first, to generate an interpolated field based on 
the covariance fields of the initial distributions and 
second, to find a probability mass function which covariance field is close to the interpolated one. 
Closeness is measured using a suitable similarity invariant function.
Covariance fields are also considered discrete ones -  
they are defined on a finite set of observation points.
With a fixed coordinate system at each observation point, not necessarily a global one, 
the covariance field is represented by a set of matrices. 
As always, we are going to use the tensor notation to guarantee a coordinate free approach.

Let $\{p_i\}_{i=1}^{k}$ and $\{q_i\}_{i=1}^{k}$ be two sets of k points on $\mathbb{S}^2$. 
The first set is the distribution domain. The second one is the observation set.
Hereafter, a discrete mass function (pmf) is any k-vector $f$, such that  
$f=\{f_i = f(p_i) \ge 0\}_{i=1}^k$ and $\sum_{i=1}^k f_i = 1$.
We write $f\in P_k^+$, where $P_k^+$ denotes the compact k-simplex. 

The number of observation points may be in fact less than k, the size of the {\it pmfs}. 
However, with a smaller observation set one may lose the uniqueness and the continuity of an estimation.
Particular geometric configurations also lead to the same result and one has  
to check carefully the consistency conditions corresponding to the problem.

The covariance field of $f\in P_k^+$ at $q_j$ is defined as 
$$
\Sigma[f]_j := \Sigma[f](q_j) = 
\sum_{i=1}^k (\overrightarrow{q_jp_i})(\overrightarrow{q_jp_i})' r(||\overrightarrow{q_jp_i}||) f(p_i),
$$
where 
$$
\overrightarrow{q_jp_i} = \frac{cos^{-1}<p_i,q_j>}{(1-<p_i,q_j>^2)^{1/2}} (p_i-<p_i,q_j>q_j).
$$
We use either $r=1$ or 
\begin{equation}\label{randcomp_eq:r_pi2}
r(t) = (1-\frac{\pi}{2t})^2.
\end{equation}
The second choice is known to be optimal on $\mathbb{S}^2$ in the class of functions 
$r_a(t) = (1-\frac{a}{t})^2$ because it minimizes the maximum of $tr(G\Sigma(q))$.

Let $f^s$,s=1,...,m, be a collection of {\it pmfs} and 
$$
\{C_j^s = \Sigma[f^s]_j\}_{j=1}^{k}, \\ s=1,...,m,
$$
be their covariance fields. 

For a non-negative similarity invariant function $h$, we define 
\begin{equation}\label{randcomp_eq:pdfs_cov_dist}
d_h(f, f^s) := \sum_{j=1}^k  h(\Sigma[f]_j, C_j^s),\\ s=1,...,m.
\end{equation}
For $\alpha\in P_m^+$, i.e. $\alpha=\{\alpha_s\}_{s=1}^m$, such that $\alpha_s \ge 0$ and $\sum_s \alpha_s = 1$, 
we define the functional
\begin{equation}\label{randcomp_eq:H}
H(f;\alpha) := \sum_{s=1}^m \alpha_s d_h(f, f^s).
\end{equation}

Then we formulate the following optimization problem: 
find a probability mass function $\hat f$ such that
\begin{equation}\label{randcomp_eq:optimal_pmf}
\hat f (\alpha) = argmin_f H(f;\alpha).
\end{equation}

Below we show some results regarding the consistency of the estimators (\ref{randcomp_eq:optimal_pmf}).

\newtheorem{randcomp_lemma1}{Lemma}
\begin{randcomp_lemma1}\label{randcomp_lemma:H_cont}
Let $h\in\mathcal{SIM}(n)$, $\alpha_s \in P_M^+$ and $f^s\in P_k^+$.
If $\alpha_s \to \alpha_0$ and $f^s \to f^0$ (in $L_2$ norm), then 
$$
H(f^s, \alpha_s) \to H(f^0, \alpha_0).
$$
\end{randcomp_lemma1}
{\it Proof.}
Observe that
$$
|H(f^s; \alpha_s) - H(f^0; \alpha_0)| \le 
|H(f^s; \alpha_s) - H(f^s; \alpha_0)| + |H(f^s; \alpha_0) - H(f^0; \alpha_0)|.
$$
Since $H(f; \alpha_0)$ is continuous in $f$, the second term above goes to zero. 
The first term is bounded by
$$
|H(f^s; \alpha_s) - H(f^s; \alpha_0)| \le ||\alpha_s - \alpha_0||_{L_2} \max_{j,m} h(\Sigma[f]_j, C_j^s).
$$
The sets $\{\Sigma[f]_j| f\in P_k^+\}$ are compact in $Sym_n^+$ and $h$ is continuous, therefore 
$\max_{j,m} h(\Sigma[f]_j, C_j^s) = C < \infty$ and 
$$
H(f^s; \alpha_s) \to H(f^s; \alpha_0).
$$
$\Box$

For a sequence $\alpha_s$, define $\hat f^s = argmin_f H(f;\alpha_s)$. 
We have the following
\newtheorem{randcomp_lemma2}[randcomp_lemma1]{Lemma}
\begin{randcomp_lemma2}\label{randcomp_lemma:Hhat_conv}
If $h\in\mathcal{SIM}(n)$ and $\alpha_s \to \alpha_0$, then 
\begin{equation}\label{randcomp_eq:H_converg}
H(\hat f^s, \alpha_s) \to H(\hat f^0, \alpha_0).
\end{equation}
\end{randcomp_lemma2}
{\it Proof.}
Since $P_k^+$ is a compact any sub sequence of $f^s$ has a point of convergence in $P_k^+$.
Without loss of generality we may assume that $\hat f^s \to g\in P_k^+$.
Accounting for the minimizing properties of $\hat f$ and applying lemma \ref{randcomp_lemma:H_cont} 
we can write 
$$
H(\hat f^0, \alpha_s) \ge H(\hat f^s, \alpha_s) \to H(g, \alpha_0) \ge H(\hat f^0, \alpha_0).
$$
Because of $H(\hat f^0, \alpha_s) \to H(\hat f^0, \alpha_0)$ we have (\ref{randcomp_eq:H_converg}).
$\Box$

Unfortunately, (\ref{randcomp_eq:H_converg}) is not enough to claim that $\hat f^s \to \hat f^0$. 
However, if $H(f;\alpha_0)$ has a well separated minimum at $\hat f^0$ we indeed have 
the wanted consistency.

\newtheorem{randcomp_corollary1}{Corollary}
\begin{randcomp_corollary1}\label{randcomp_lemma:fhat_cont}
$\hat f(\alpha)$ is continuous at all $\alpha$ for which $H(f, \alpha)$ has a well separated (global) minimum.
\end{randcomp_corollary1}

Another problem is how to find the global minimum $\hat f$ of $H(f;\alpha)$, provided it is unique.
We know that the minimum is easily found in case of convex function H, by gradient descent algorithm for example.
Moreover, the convexity of $H(f;\alpha_0)$ in $P_k^+$ guarantees the well separability of its minimum and 
that gives us the desired consistency.

\newtheorem{randcomp_prop1}{Proposition}
\begin{randcomp_prop1}\label{randcomp_prop:f_conv}
If $\alpha_s \to \alpha_0$ and $h\in \mathcal{SIM}(n)$ is such that $H(f;\alpha_0)$ 
is convex in $P_k^+$, 
then 
\begin{equation}\label{randcomp_eq:f_converg}
\hat f^s \to \hat f^0.
\end{equation}
\end{randcomp_prop1}
{\it Proof.}
Suppose the contrary of (\ref{randcomp_eq:f_converg}), that there exists $g\in P_k^+$, and 
sub sequence $\hat f^s \to g$, such that $||\hat f^0 - g|| > 0$. 
Then $H(g;\alpha_0) > H(\hat f^0; \alpha_0)$ by the separability of the minimum. 
But  
$H(\hat f^s; \alpha_s) \to H(g; \alpha_0)$ by lemma \ref{randcomp_lemma:H_cont} 
and
$H(\hat f^s; \alpha_s) \to H(\hat f^0; \alpha_0)$ by lemma \ref{randcomp_lemma:Hhat_conv}, 
which imply $H(g;\alpha_0) = H(\hat f^0; \alpha_0)$. The contradiction proves the claim.
$\Box$

%%%%%%%%%%%%%%%%%%%%%%%%%%%%%%%%%%%%%%%%%%%%%%%%%%
{\subsection{Linear Interpolation}}

Consider first one of the simplest similarity invariant functions $h_{trdif}^2(.,.; G^{-1})$. 
The corresponding optimization functional is
$$
H_{trdif}(f, \alpha) = \sum_{s=1}^m \alpha_s \sum_{j=1}^{k} tr^2(G(q_j)\Sigma[f]_j - G(q_j)C_j^s)
$$
Denote $a_{ij} = tr(G(q_j)(\overrightarrow{q_jp_i})(\overrightarrow{q_jp_i})') = d^2(q_j, p_i)$ and $c_j^s = tr(G(q_j)C_j^s)$, 
then 
$$
H_{trdif}(f, \alpha) = \sum_{s=1}^m \alpha_s \sum_{j=1}^{k} (\sum_i a_{ij}f_i - c_j^s)^2.
$$
We have
$$
\frac{\partial H_{trdif}}{\partial f_i} = 2\sum_{s=1}^m \alpha_s \sum_{j=1}^{k} a_{ij}(\sum_l a_{lj}f_l - c_j^s).
$$
The second partial derivatives are 
$$
\frac{\partial^2 H_{trdif}}{\partial f_i\partial f_l} = 2\sum_{s=1}^m \alpha_s \sum_{j=1}^{k} a_{ij} a_{lj}
$$
Let $w=\{w_i\}\in \mathbb{R}^k$, then
$$
\sum_{i,l} w_i w_l \frac{\partial^2 H_{trdif}}{\partial f_i\partial f_l} = 
2\sum_{s=1}^m \alpha_s \sum_{j=1}^k (\sum_{i=1}^k w_i a_{ij}) ^ 2 \ge 0.
$$
Therefore, if the matrix $A = \{a_{ij}\}_{i=1,j=1}^{k,k}$ is of full rank $k$, then 
$H_{trdif}$ is convex in $P_k^+$.
Moreover, the optimal solution of (\ref{randcomp_eq:optimal_pmf}) satisfies 
$$
\sum_i a_{ij}f_i = \sum_{s=1}^m \alpha_s c_j^s,\\ j=1,...,k,
$$
with a unique solution 
$$
\hat f = \sum_{s=1}^m \alpha_s f^s,
$$
since for every $s$ and $j$, $\sum_i a_{ij}f_i^s = c_j^s$.

Thus, we showed the following
\newtheorem{randcomp_prop2}[randcomp_prop1]{Proposition}
\begin{randcomp_prop2}\label{randcomp_prop:linear_interpolation}
If the matrix A has full rank, $rank(A) = k$, then the linear interpolation is the unique solution 
of the optimization problem (\ref{randcomp_eq:optimal_pmf}) 
for $H_{trdif}$.
\end{randcomp_prop2}

%%%%%%%%%%%%%%%%%%%%%%%%%%%%%%%%%%%%%%%%%%%%%%%%%%
{\subsection{Non-Linear Interpolations}}

Consider similarity invariant function $h_{trln2}$ and corresponding optimization functional $H_{trln2}$
$$
H_{trln2}(f; \alpha) = 
\sum_{s=1}^m \alpha_s \sum_{j=1}^{k} tr(ln^2(\Sigma[f]_j(C_j^s)^{-1})).
$$
The value of $H_{trln2}(f)$ is small when $G(q_j)\Sigma[f]_j$ is close to 
covariance operators $G(q_j)C_j^s$ for all $j$ and $s$.
This is a much stronger condition than the requirement for their traces to be close as in the problem of 
minimizing $H_{trdif}$. Consequently the minimum of $H_{trln2}(f)$, in general, will be strictly positive and 
the optimal $pmf$ will be different from the linear interpolation.

Experiments show great improvement in convergence of gradient descend algorithm 
for problem (\ref{randcomp_eq:optimal_pmf}), when instead of the generic covariance one uses 
the second choice (\ref{randcomp_eq:r_pi2}).

Define the operators 
$$
Z_{ij}^s = (\overrightarrow{q_jp_i})(\overrightarrow{q_jp_i})'(1-\frac{\pi}{2||\overrightarrow{q_jp_i}||})^2 (C_j^s)^{-1}
$$
and set $Y_j^s = \sum_i f_i Z_{ij}^s$. 
The gradient of $H_{trln2}$ is 
$$
\nabla H_{trln2}(f, \alpha) = 
\{ f_i \sum_{s=1}^m \alpha_s \sum_{j=1}^{k} \frac{tr(ln(Y_j^s)Z_{ij}^s)}{tr(Z_{ij}^s)} \}_{i=1}^k.
$$
The optimization problem (\ref{randcomp_eq:optimal_pmf}) is solved by gradient descent algorithm, 
which shows relatively fast convergence,
unfortunately not always to the global minimum, because $H_{trln2}(f, \alpha)$ is not convex in $f\in P_k^+$.

Log-likelihood function gives us another choice for $H$, 
$$
H_{lik}(f; \alpha) = 
\sum_{s=1}^m \alpha_s \sum_{j=1}^{k} \{tr(\Sigma[f]_j(C_j^s)^{-1}) - ln|\Sigma[f]_j(C_j^s)^{-1}| - n\} = 
$$
$$
\sum_{s=1}^m \alpha_s \sum_{j=1}^{k} \{tr(Y_j^s) - ln|Y_j^s| - n\}.
$$
The gradient of $H_{lik}$ is  
$$
\nabla H_{lik}(f; \alpha) = 
\{ f_i \sum_{s=1}^m \alpha_s \sum_{j=1}^{k} \frac{tr((Y_j^s-I_n)Z_{ij}^s)}{tr(Z_{ij}^s)} \}_{i=1}^k.
$$

Note that $h_{lik}$ is neither symmetric nor satisfies the triangular inequality, 
but its importance is determined by the relation to normal distributions and its analytical properties. 
Define the matrix 
$$
B = \{ b_{ij} = (d(q_j,p_i)-\frac{\pi}{2})^2 \}_{i=1,j=1}^{k,k}.
$$
\newtheorem{randcomp_prop3}[randcomp_prop1]{Proposition}
\begin{randcomp_prop3}\label{randcomp_prop:loglik_convex}
If B has full rank, $rank(B)=k$, then for all $\alpha$, $H_{lik}(f;\alpha)$ is a convex function in $P_k^+$.
\end{randcomp_prop3}
{\it Proof.}
We have
$$
\frac{\partial H_{lik}}{\partial f_i} = \sum_{s=1}^m \alpha_s \sum_{j=1}^{k} 
tr(Z_{ij}^s - Z_{ij}^s (Y_j^s)^{-1}).
$$
and
$$
\frac{\partial^2 H_{lik}}{\partial f_i\partial f_l} = \sum_{s=1}^m \alpha_s \sum_{j=1}^{k} 
tr(Z_{ij}^s(Y_j^s)^{-1}Z_{lj}^s(Y_j^s)^{-1}).
$$
We want to show that the matrix of second partial derivatives is positive definite.
Let $w=\{w_i\}\in \mathbb{R}^k$ and $w\ne 0$, then
$$
\sum_{i,l} w_i w_l \frac{\partial^2 H_{lik}}{\partial f_i\partial f_l} = 
\sum_{s=1}^m \alpha_s \sum_{j=1}^k tr(\sum_{i=1}^k w_i Z_{ij}^s(Y_j^s)^{-1}) ^ 2 > 0,
$$
since by the assumption for $B$, for at least one $j$, $\sum_{i=1}^k w_i Z_{ij}^s \ne 0$.
$\Box$

The rank of $B$ can be calculated using the pairwise distances between $q$ and $p$ points 
and only in very special circumstances this rank will be less than $k$. 
More formally, if a random process chooses the points, then 
$$
P(rank(B) < k) = 0.
$$

%%%%%%%%%%%%%%%%%%%%%%%%%%%%%%%%%%%%%%%%%%%%%%%%%%
{\subsection{Examples and conclusions}}

\begin{figure}
\centering
\begin{tabular}{cc}
{\includegraphics[scale=0.4]{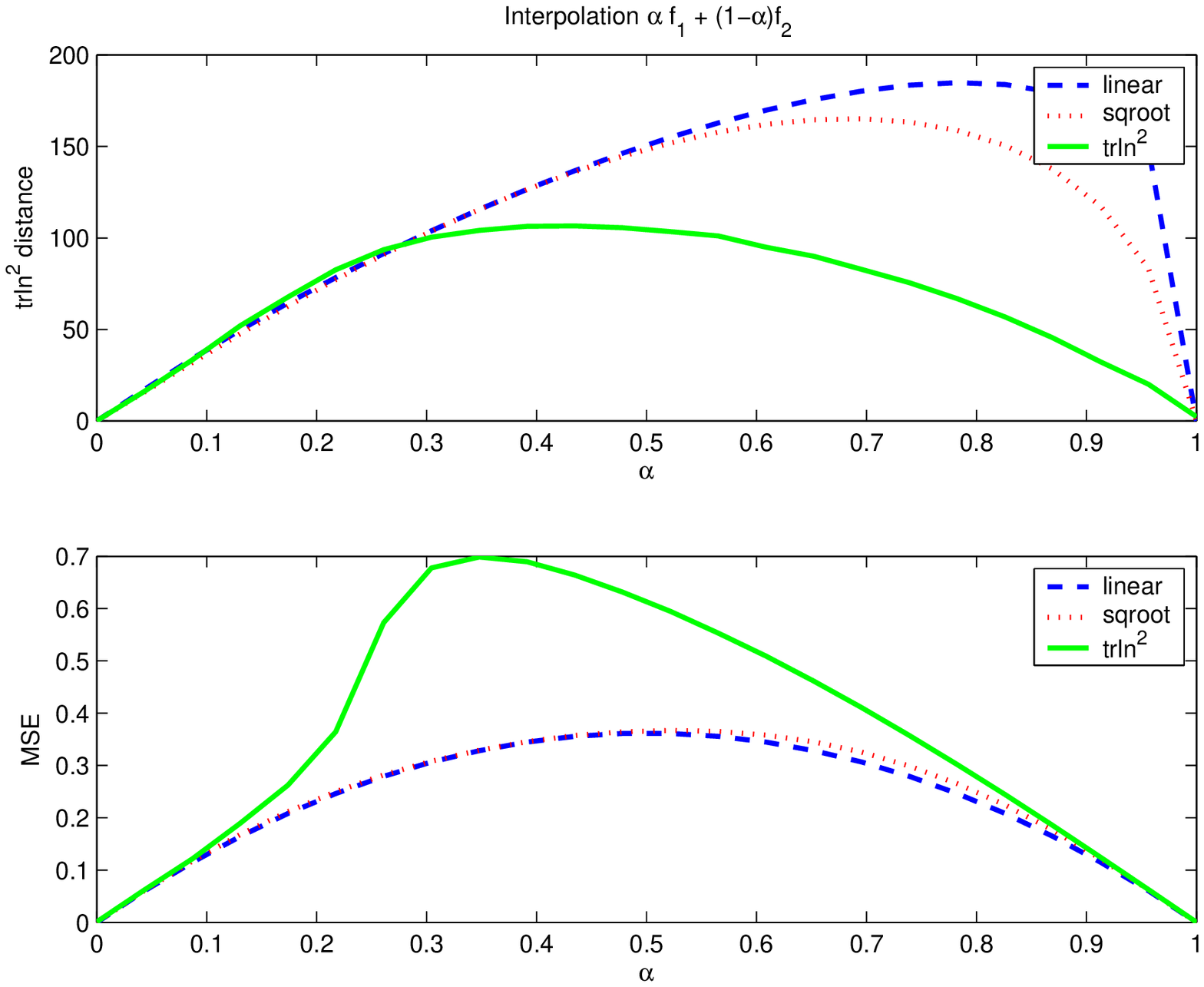}}
{\includegraphics[scale=0.4]{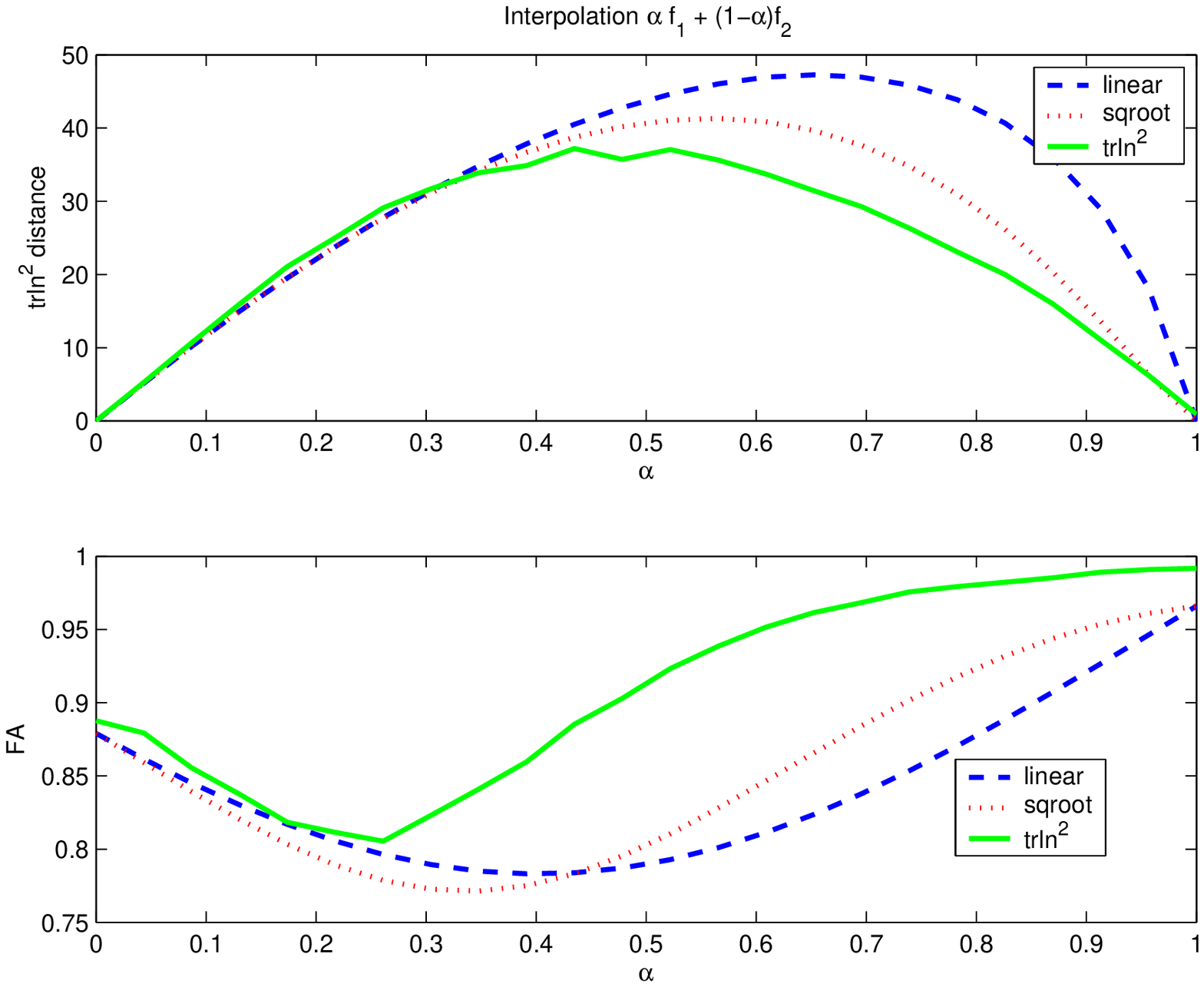}}
\end{tabular}
\caption{Two examples of interpolation of {\it pmfs} on $\mathbb{S}^2$ using $h_{trln2}$. 
The linear and square root (see (\ref{randcomp_eq:sqroot}))
interpolations are also given for reference. 
Top plots show $H_{trln2}$ and $H_{lik}$ for the three interpolations.
Bottom plots show corresponding MSEs (see (\ref{randcomp_eq:MSE})) in the left and 
FAs (see (\ref{randcomp_eq:FA})) in the right. 
}
\label{randcomp_fig:7}
\end{figure}

Figure \ref{randcomp_fig:7} shows interpolation between two $pmfs$ of size 6 ($m=2, k=6$) 
applying $h_{trln2}$. We compare it to the linear and the square root interpolations. 
Square root interpolation, as suggested by the name, relies on the observation that 
for a $pmf$ $f\in P_k^+$, $\sqrt{f}=(\sqrt{f^1},...,\sqrt{f^k}) \in \mathbb{S}^k$. Then 
one finds 
\begin{equation}\label{randcomp_eq:sqroot}
\hat p = argmin_{p\in \mathbb{S}^k} \sum_s\alpha_s d^2(p,\sqrt{f^s}) \textrm{ and sets } \hat f_{sqroot} = \hat{p}^2.
\end{equation}
It is also informative to compare the Mean-Squared Error (MSE) between different interpolations. 
It is defined by
\begin{equation}\label{randcomp_eq:MSE}
MSE(\hat f) = \sum_{s=1}^2 \alpha_s \sum_{i=1}^k (\hat f_i - f_i^s)^2.
\end{equation}
Linear and square root interpolations, by their nature, are very close in MSE, 
but very different from $\hat f_{trln^2}(\alpha)$, which manifests the non-linear origin of the latter.

Another performance criteria relevant to the study of spherical data is the Fractional Anisotropy (FA).
Let $\{\lambda_i\}_{i=1}^n$ be the eigenvalues of $\sum_{i=1}^k \overrightarrow{p_i}\overrightarrow{p_i}'f_i$, 
where $\overrightarrow{p_i}$ are considered vectors in $\mathbb{R}^3$ 
(thus $FA$ is defined only for distributions on $\mathbb{S}^2$).
Then we define
\begin{equation}\label{randcomp_eq:FA}
FA(f) = \{\frac{n}{n-1}\sum_{i=1}^n(\lambda_i-\bar\lambda)^2 / \sum_{i=1}^n\lambda_i^2 \}^{1/2}.
\end{equation}
Fractional Anisotropy measures a distribution concentration. The higher FA the more concentrated it is about particular axes.
A uniform distribution has $FA=0$. 
As we may expect the linear interpolation substantially reduces the FA index, 
$h_{trln2}$-based one however, is more conservative and manage to sustain 
higher FA. Preserving the concentration factor is of importance for processing ODFs in HARDI, 
and the empirical evidence for the good FA performance of $h_{trln2}$ is encouraging.

\begin{figure}
\centering
\begin{tabular}{cc}
{\includegraphics[scale=0.4]{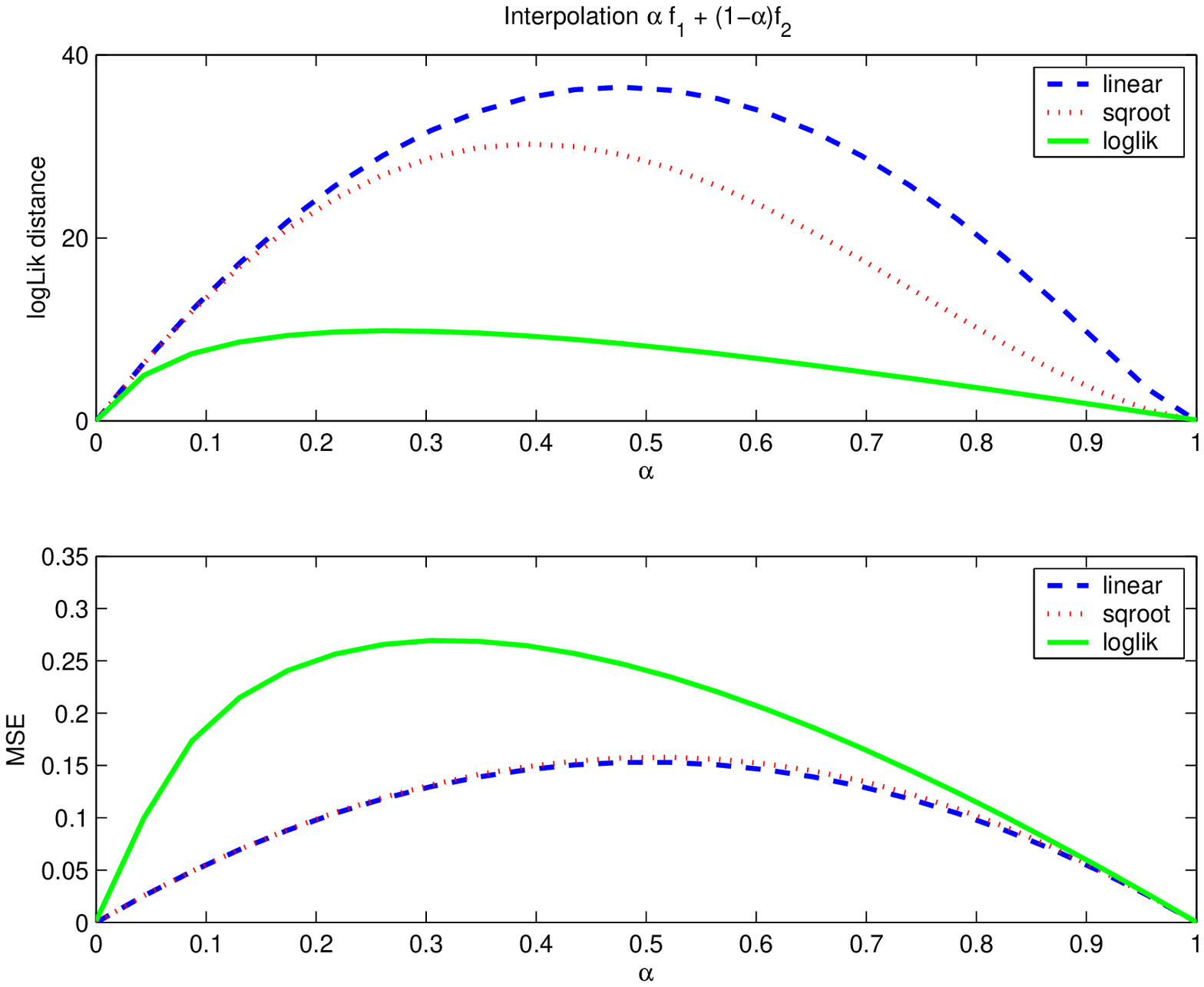}}
{\includegraphics[scale=0.4]{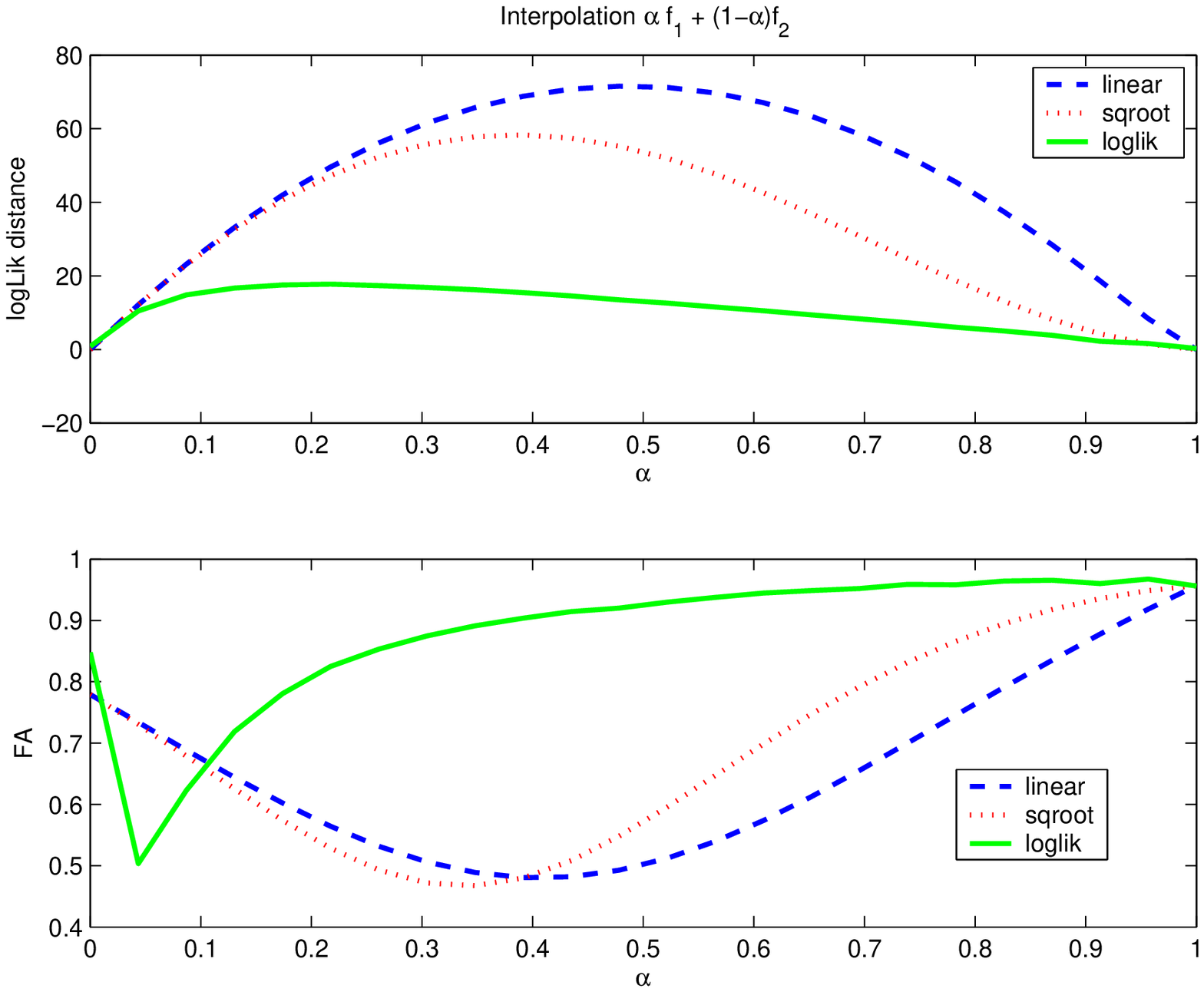}}
\end{tabular}
\caption{Two examples of interpolation of {\it pmfs} on $\mathbb{S}^2$ using $h_{lik}$. 
The linear and square root 
interpolations are also given for reference. 
Top plots show $H_{lik}$ for the three alternatives.
Bottom plots show corresponding MSEs in the left and 
FAs in the right.
}
\label{randcomp_fig:8}
\end{figure}

A second set of examples in figure \ref{randcomp_fig:8} illustrates interpolation based on 
the likelihood function, $h_{lik}$. As we showed, this choice guarantees the convexity of 
$H_{lik}$ and thus the continuity of the optimal solution $\hat f_{lik}(\alpha)$.

The likelihood based interpolation $\hat f_{lik}$ exhibits behaviour similar to that of $\hat f_{trln2}$. 
Again, it is very distinguished from the linear and the square-root one and tends to preserve the anisotropy.

%

%%%%%%%%%%%%%%%%%%%%%%%%%%%%%%%%%%%%%%%%%%%%%%%%%%
{\section{Summary}}

The main goal of this article is to introduce covariance operator fields and provide 
some arguments showing their potential and usefulness.

There is a covariance field associated with any distribution on a Riemannian manifold. 
It defines a linear operator on the tangent space of each point on the manifold. 
By applying a similarity invariant to that operator field one can obtain a scalar field 
that represents the distribution. 
It reveals important spatial characteristics of the distribution.
Similarity invariants can also be used for comparing and interpolating distributions.

We demonstrated several non-parametric procedures for solving 
a two-sample location problem on the sphere and showed how covariance operator fields 
can be used for locating observation points that maximize test performance.

We also implemented two non-linear procedures for interpolating distributions 
on the sphere and compared them to the linear and square-root interpolations. 
The proposed approach is general enough to allow a great variety of choices and promises 
a good application potential. 

%%%%%%%%%%%%%%%%%%%%%%%%%%%%%%%%%%%%%%%%%%%%%%%%%%
{\section{Acknowledgements}}

We thank Vic Patrangenaru for the helpful comments and useful discussions.
We also thank Paul Thompson from UCLA School of Medicine for providing us with a HARDI volume for research purposes.

\end{document}